\documentclass[12pt,thmsa]{amsart}

\usepackage{amssymb, euscript}
\def\Cal{\mathcal}

\def\P{{\Cal P}}

\def\F{{\Cal F}}

\def\I{{\Cal I}}
\def\L{{\Cal L}}

\def\W {{\EuScript W}}

\def\S{\EuScript{S}}
\def\nf{\|f\|_p}

\def\tr{{\hbox{\rm tr}}}

\def\Ma{\frM_{n,m}}

\def\Mt{\frM_{n-k,m}}
\def\Mkm{\frM_{k,m}}

\def\nf{\|f\|_p}
\def\V{{\bf W}_{n,m}}

\def\f0{f_0}
\def\Fc0{\varphi_0}

\def\irn{\intl_{\bbr^n}}

\def\I_k {I_{-}^{k/2}}
\def\I+k {I_{+}^{k/2}}

\def\vnk{V_{n,n-k}}
\def\vnm{V_{n,m}}
\def\cd{\vnk\times\frM_{n-k,m}}
\def\Gr{\frT}

\def\bbr{{\Bbb R}}

\def\bbc{{\Bbb C}}

\def\rank{{\hbox{\rm rank}}}

\def\tr{{\hbox{\rm tr}}}

\def\const{{\hbox{\rm const}}}

\def\det{{\hbox{\rm det}}}

\def\min{{\hbox{\rm min}}}

\def\fc{\varphi(\xi,t)}

\def\lp{L^p(\frM_{n,m})}

\def\part{\partial}
\def\intl{\int\limits}
\def\b{\beta}

\def\Gam{\Gamma}

\def\a{\alpha}
\def\om{\omega}

\def\del{\delta}
\def\vp{\varphi}

\def\g{\gamma}
\def\gam{\gamma}

\def\sig{\sigma}
\def\lam{\lambda}
\def\z{\zeta}
\def\e{\varepsilon}
\def\t{\tau}
\def\eq{\xi 'x=t}

\def\rf{\hat f (\xi,t)}
\def\df{\check \varphi (x)}

\font\frak=eufm10

\def\fr#1{\hbox{\frak #1}}

\def\frM{\fr{M}}

\def\frT{\fr{T}}

\def\const{{\hbox{\rm const}}}

\def\det{{\hbox{\rm det}}}

\def\min{{\hbox{\rm min}}}
\def\span{{\hbox{\rm span}}}

\def\p{\P_m}
\def\gm{\Gamma_m}
\def\tr{{\hbox{\rm tr}}}
\def\part{\partial}
\def\intl{\int\limits}
\def\b{\beta}

\def\Gam{\Gamma}

\def\a{\alpha}
\def\cpm{\overline\P_m}
\def\pd{\stackrel{*}{P}\!{}^\a}

\newtheorem{theorem}{Theorem}[section]
\newtheorem{lemma}[theorem]{Lemma}

\theoremstyle{definition}
\newtheorem{definition}[theorem]{Definition}

\newtheorem{corollary}[theorem]{Corollary}
\theoremstyle{remark}

\usepackage{amssymb}

 \numberwithin{equation}{section}

\newcommand{\be}{\begin{equation}}
\newcommand{\ee}{\end{equation}}

\newcommand{\bea}{\begin{eqnarray}}
\newcommand{\eea}{\end{eqnarray}}
\newcommand{\Bea}{\begin{eqnarray*}}
\newcommand{\Eea}{\end{eqnarray*}}

\begin{document}

\title [Multiscaled  wavelet transforms  ]{ Multiscaled  wavelet transforms, ridgelet
transforms, and  Radon transforms on the space of matrices}

\author{G. \'Olafsson}

\address{Department of Mathematics, Louisiana State University, Baton
Rouge, Louisiana 70803, USA} \email{olafsson@math.lsu.edu}

\author{ E. Ournycheva}
\address{Institute of Mathematics,
Hebrew University, Jerusalem 91904, \newline ISRAEL}
\email{ournyce@math.huji.ac.il}

\author{ B. Rubin}
\address{Institute of Mathematics,
Hebrew University, Jerusalem 91904, \newline ISRAEL
\newline
and
\newline
Baton Rouge, LA 70803,
USA}
\email{boris@math.huji.ac.il}
\thanks{GO was supported
by NSF Grants DMS-0402068 and DMS-0139783}
\thanks{BR was supported by the Edmund Landau Center for Research in Mathematical Analysis
and Related Areas, sponsored by the Minerva Foundation (Germany)
and  partially by DMS-0139783}

\subjclass[2000]{Primary 42C40; Secondary 44A12}

\date{September 2, 2004.}


\keywords{The Radon  transform, matrix spaces, the  Fourier
transform, Riesz potentials, wavelet transforms, ridgelet
transforms}

\begin{abstract}
 Let  $\Ma$ be the space
of real $n\times m$ matrices which can be identified with   the
Euclidean space $\bbr^{nm}$.  We introduce continuous wavelet
transforms on $\Ma$ with a multivalued  scaling parameter
represented by a positive definite symmetric  matrix. These 
transforms agree with the  polar decomposition on $\Ma$ and
coincide with  classical ones in the rank-one case $m=1$. We prove
an analog of Calder\'{o}n's reproducing formula for
$L^2$-functions   and obtain explicit  inversion
formulas for the Riesz potentials and Radon
transforms on  $\Ma$. We also introduce continuous ridgelet transforms
associated to matrix planes in $\Ma$. An inversion formula for
these transforms follows from that for the Radon transform.
\end{abstract}

\maketitle

\section{Introduction}

\setcounter{equation}{0}

\noindent It is known that diverse wavelet-like transforms can be generated
by operators of fractional integration and used to  invert these operators. 
  On the other hand, numerous problems in 
integral geometry, for instance, reconstruction of functions from
their integrals over planes in $\bbr^n$, reduce to inversion of  fractional integrals. The following  example illustrates these statements and explains 
how  wavelet transforms arise  in the context of
integral-geometrical problems; see  also \cite{Ru2} for the more detailed exposition.

 Consider the
Riesz potential \be\label {ri} (I^\a f)(x) \, = \, {1\over
\gamma_n (\a)} \intl_{\bbr^n} |x-y|^{\a-n}f(y) dy, \qquad x \in
\bbr^n, \ee
\[  \gamma_n (\a) = \; {\pi^{n/2} 2^\a \Gamma (\a/2)\over
\Gamma ((n-\a)/2)},  \qquad Re \, \a >0, \quad \a -n\neq 0,2,
\dots \, ,
\]
and replace the kernel $ |x-y|^{\a-n}$ by the integral
\be\label{ker} |x-y|^{\a-n} = c^{-1}_{\a, w} \intl^\infty_0 w
\left({|x-y|\over a}\right) \; {da\over a^{n-\a+1}}, \ee
where  $w(\cdot)$ is good enough and 
$c_{\a, w} = \int^\infty_0 w(s)\; s^{n-\a-1}
 ds \neq 0$.
 Changing the
order of integration, we obtain \be\label{pot} (I^\a f)(x) \, =\,
{c^{-1}_{\a, w}\over \gamma_n(\a)} \, \intl^\infty_0 \, {(\W_a
f)(x)\over a^{1-\a}}\, da, \ee where
\be \label{wav}(\W_af)(x) = \frac{1}{a^n} \intl_{\bbr^n} f(y) w
\left(\frac{|x-y|}{a}\right) \; dy. \ee If  $w$ obeys some
cancellation conditions, then  (\ref{wav}) represents the classical
continuous wavelet transform with the scaling parameter $a>0$
\cite{Da}, \cite{FJW}, \cite {Ho}. If we start  with a fractional integral 
different from $I^\a f$, say, with the Bessel potential or
 whatever (see \cite{Ru1}, Section 10.7), we arrive at  a wavelet transform, 
 which differs from (\ref{wav}).

Since the Fourier transform of $I^\a f$ is $|y|^{-\a} (\F f)(y)$
in a certain sense, then,   formally, $(I^\a)^{-1} =I^{-\a}$, and
it is natural to expect that  the inverse operator
$(I^\a)^{-1}$ can be represented in the form (\ref{pot}) with $\a$
replaced by $-\a$, namely, 
\be \label{inv} (I^\a)^{-1} f =
d_{\a, w} \intl^\infty_0 \frac{\W_af}{ a^{1+\a}} \, da, \ee
$d_{\a, w}$ being a normalizing factor. In particular, for $\a=0$,
\be\label{cal} f = d_{0, w} \intl^\infty_0 \frac{\W_af}{ a} \, da. \ee The equality (\ref{cal}) is a modification of Calder\'on's
reproducing formula \cite{Cal}, \cite{FJW}, \cite{Ru3},
\cite{Ru7}. If $w$ is normalized and has the form $w=u \ast v$,
then (\ref{cal}) turns into  the classical Calder\'on identity
\be\label {clas} f = \intl^\infty_0 \frac{f \ast u_a \ast
v_a}{a} \, da \ee where $u_a(x)=a^{-n }u(x/a), \; v_a(x)=a^{-n
}v(x/a)$.

Of course, this argument is purely heuristic and formulas
(\ref{inv})-(\ref{clas}) require justification in the framework of
a suitable class of functions $f$ under certain cancellation
conditions for the wavelet function $w$.

What is the  connection between this argument and the Radon
transform in integral geometry? Suppose that $\frT$ is the
manifold of all $k$-dimensional planes $\t$ in $\bbr^n$. For functions
$f : \bbr^n \to \bbc$ and $\vp : \frT \to \bbc$, the Radon
transform and its dual are defined by \be \label{rd} \hat f
(\t)=\intl_{x \in \t} f(x), \qquad \check \vp (x)=\intl_{\t \ni x}
\vp(\t), \ee respectively, and obey the Fuglede equality
\be\label{fu} (\hat f)^\vee =c \, I^k f, \qquad c=\const,\ee see
\cite {Fu}, \cite{Hel}, \cite{Ru4}  for details. Combining
(\ref{fu}) with (\ref{inv}), we obtain an inversion formula for
$\hat f$ in the ``wavelet form" \be \label {irw} f=\const
\intl^\infty_0 \frac{\W_a (\hat f)^\vee }{a^{1+k}} \, da, \ee
where $\int_0^\infty=\lim_{\e \to 0} \int_\e^\infty$ in a certain
sense.

This formalism can be applied in a wider context, when, instead of
the Euclidean distance $|x-y|$ between two points, one deals with
the distance $|x-\tau|$ between the point $x \in \bbr^n$ and the
$k$-dimensional plane $\tau \subset \bbr^n$. Starting with the
intertwining operator  \be \label{sem}(P^\a
f)(\t)=\frac{1}{\gamma_{n-k}(\a)}\irn f(x) |x-\t|^{\a+k-n} \, dx
\ee which is  called  the {\it  generalized Semyanistyi fractional
integral} (cf. \cite {Se} for $k=n-1$),
 one arrives at the
corresponding wavelet-like transform, recently called the {\it
continuous $k$-plane ridgelet transform}; see   \cite{Ca}, \cite{Ru5},  and
references therein. These transforms have proved to be useful in
applications \cite{Ca}, \cite{Do}, \cite {Mur}, and are of
independent theoretical interest.

A common feature of these examples is that the scaling parameter
$a>0$ is one-dimensional  no matter what the dimension of the
ambient space $\bbr^n$ is.  The situation changes drastically  if
we replace $n$ by
 $nm$, and regard $\bbr^{nm}$ as the space
  of $n \times m$
real matrices  $x=(x_{i,j})$. Then a similar procedure, starting with the properly defined 
 Riesz potential, yields a new wavelet-like
transform which is applicable to  functions of matrix argument and
depends on the matrix-valued scaling parameter. This parameter is
represented by a positive definite symmetric matrix of size $m
\times m$. Apart from extra flexibility that might be useful in
applications, such  transforms have a rich theory which relies
on diverse higher rank phenomena.

In the present paper, we focus on the $L^2$ theory
of the  new wavelet
transforms mentioned above, and give some applications. The paper is organized as
follows.  Section 2 contains necessary prerequisites. We fix 
our notation and recall basic
facts related to Riesz potentials and Radon transforms on the
space of rectangular matrices. In Section 3, we introduce
continuous wavelet transforms for functions of matrix argument
and prove the corresponding reproducing formula of the
Calder\'on type. In Section 4, we show how  wavelet transforms
can be used for inversion of Riesz potentials on matrix spaces.
Unlike the rank-one case $m=1$, for which numerous inversion
formulas are known \cite{Ru1}, \cite{SKM}, the corresponding higher rank problem is
very difficult; see \cite{OR}, \cite{Ru6}  for the discussion. Wavelet transforms
prove to be a convenient tool to resolve this problem in the $L^2$-case. In Section 5, we
apply our wavelet transforms to
inversion of the Radon transform associated to the so-called
matrix $k$-planes. These  Radon transforms were studied in detail in
\cite{Pe},  \cite{OR1}, and  \cite{OR}, where it was shown that the inversion problem for
them has the same difficulties as for the Riesz potentials on the space of 
rectangular matrices. In Section 5, we
introduce continuous ridgelet transforms of  functions of matrix
argument, generalizing those in  \cite{Ca} and \cite{Ru5}, and prove a
reproducing formula for these transforms. This result is a
consequence of the inversion formula for the   Radon
transform.
\medskip

\noindent
{\bf Acknowledgements.} The work on the paper was started when B. Rubin was visiting
Department of Mathematics at the Louisiana State University in April 2004. He
is deeply
grateful to his colleague, Prof. Gestur \'Olafsson for the hospitality.

\section{Preliminaries}

\setcounter{equation}{0}

\noindent In this section, we fix  our notation and recall some basic facts,
 that will be used
throughout the paper. The main references are \cite{Mu},
\cite{OR}, \cite{T}.
\subsection{Notation and some auxiliary facts}
  Let $\frM_{n,m}$ be the
space of real matrices $x=(x_{i,j})$ having $n$ rows and $m$
 columns. We identify $\frM_{n,m}$
 with the real Euclidean space $\bbr^{nm}$ and set $dx=\prod^{n}_{i=1}\prod^{m}_{j=1}
 dx_{i,j}$ for the Lebesgue measure on $\Ma$. In the following,
  $x'$ denotes the transpose of  $x$, $I_m$ is the identity $m \times m$
  matrix,  $0$ stands for zero entries. Given a square matrix $a$,  we 
denote by $\tr (a)$  the trace of $a$, and by $|a|$ the absolute value of
 the determinant of $a$,  respectively. We hope the reader will not confuse
 $|a|$ with the similar notation for the absolute value of a number because 
the meaning of $a$ will be clear each time from the context. 
For $x \in \Ma$, $n\geq m$, we set
 \be
 |x|_m =\det (x'x)^{1/2}= |x'x|^{1/2}.\ee If $m=1$, this is the usual Euclidean norm in
$\bbr^n$. If $m>1$, then  $|x|_m$ is the volume of the
parallelepiped spanned  by the column-vectors of the matrix $x$, cf. 
\cite[ p. 251]{G}.

 Let  $\p$ be the cone
of positive definite symmetric matrices $r=(r_{i,j})_{m\times m}$
with the elementary volume $ dr=\prod_{i \le j} dr_{i,j}$, and let
$\cpm$ be  the closure of $\p$, that is the set of  all positive
semi-definite $m\times m$ matrices. We write  $r>0$ if 
 $r\in\p$, and $r\geq 0$
 if  $r\in\cpm$, respectively. Given $s_1$ and  $s_2$ in $\cpm$, we write
 $s_1 > s_2$  for  $s_1 - s_2 \in
 \p$.
  If $a\in\cpm$ and $b\in\p$, then   $\int_a^b f(s) ds$ denotes
 the integral over the compact set
$$
(a,b)= \{s : s-a \in \p, \, b-s \in \p \}.
$$
The group $G=GL(m,\bbr)$ of
 real non-singular $m \times m$ matrices $g$ acts transitively on $\p$
  by the rule $r \to grg'$.  The corresponding $G$-invariant
 measure is  \be\label{2.1}
  d_{*} r = |r|^{-d} dr, \qquad |r|=\det (r), \qquad d= (m+1)/2 \ee \cite[p.
  18]{T}.

A function $w_0$ on $\p$ is called {\it symmetric } if
\be\label{sym} w_0( s^{1/2}rs^{1/2})=w_0(r^{1/2}sr^{1/2}),\qquad
\forall r,s\in\p. \ee For example, any function of the form
$w_0(r)=w_1(\tr(r))$ is symmetric. Instead of the trace, one can
take any function of the form $\psi(\sig_1,\dots,\sig_m)$, where
$\sig_1,\dots,\sig_m$ are elementary symmetric functions of the 
eigenvalues $\lam_1,\dots, \lam_m$ of $r$. This follows from the
general fact that if $A$ and $B$ are nonsingular square matrices
then $AB$ and $BA$ have the same eigenvalues; see, e.g., \cite[pp.
584, 585]{Mu}.

We use a standard notation $O(n)$   for the group of
real orthogonal $n\times n$ matrices; $SO(n)=\{\gam  \in O(n): 
\det(\gam)=1\}$. The corresponding invariant
measures on $O(n)$ and $SO(n)$ are normalized to be of total mass
1. The  Lebesgue space $L^p=\lp$
 and the Schwartz space  $\S=\S(\Ma)$ are
 identified with  respective spaces on $\bbr^{nm}$. We denote by $C_c(\Ma)$  the space of compactly supported
continuous functions on $\Ma$.

The Fourier transform  of a
function $f\in L^1(\Ma)$ is defined by \be\label{ft} (\F
f)(y)=\intl_{\Ma} \exp(\tr(iy'x)) f (x) dx,\qquad y\in\Ma \; .\ee
This is the usual Fourier transform on $\bbr^{nm}$ so that the
relevant Parseval formula  reads 
\be\label{pars} (\F f, \F
\vp)=(2\pi)^{nm} \, (f,\vp),\ee where
$$(f, \vp)=\intl_{\Ma} f(x) \, \overline{ \vp(x)} \, dx.$$

 We write $c$, $c_1$, $c_2, \dots$ for different constants
the meaning of which is clear from the context.
\begin{lemma}\label{12.2} [\rm (see, e.g.,
 [Mu, pp. 57--59])]
 \hfill
 
 \noindent
 {\rm (i)} \ If $ \; x=ayb$, where $y\in\Ma, \; a\in  GL(n,\bbr)$, and $ b \in  GL(m,\bbr)$, then
 $dx=|a|^m |b|^ndy.$\\
 {\rm (ii)} \ If $ \; r=qsq'$, where $s\in \P_m$ and $ q\in  GL(m,\bbr)$,
  then $dr=|q|^{m+1}ds.$ \\
  {\rm (iii)} \ If $ \; r=s^{-1}$,  $s\in \p$,   then $r\in
  \p$
  and $dr=|s|^{-m-1}ds.$
\end{lemma}

 The   Siegel gamma function
 associated to the cone $\p$ is defined by
\be\label{2.4}
 \gm (\a)=\intl_{\p} \exp(-\tr (r)) |r|^{\a -d} dr, \qquad d=(m+1)/2,\ee
 \cite{Si},   \cite{Gi}, \cite{Mu}, \cite{FK}, \cite{T}.
This integral  converges absolutely
 if and only if $Re \, \a>d-1$, and can be written as a product of
 ordinary $\Gamma$-functions:
\be\label{2.5}
 \gm (\a)=\pi^{m(m-1)/4} \Gam (\a) \Gam (\a- \frac
{1}{2}) \ldots \Gam (\a- \frac {m-1}{2}). \ee

For $n\geq m$, let $\vnm= \{v \in \frM_{n,m}: v'v=I_m \}$
 be the Stiefel manifold of orthonormal $m$-frames in $\bbr^n$.
 We fix the invariant measure $dv$ on
 $\vnm$ [Mu, p. 70]
  normalized by \be\label{2.16} \sigma_{n,m}
 \equiv \intl_{\vnm} dv = \frac {2^m \pi^{nm/2}} {\gm
 (n/2)}\; \ee and denote  $d_\ast v=\sig^{-1}_{n,m} dv$.
A polar decomposition on $\Ma$ is defined according to the
following lemma; see, e.g.,  \cite[pp. 66, 591]{Mu}, \cite{Ma}.
\begin{lemma}\label{l2.3} Let $x \in \frM_{n,m}, \; n \ge m$. If  $\rank (x)= m$,
then
\[ x=vr^{1/2}, \qquad v \in \vnm,   \qquad r=x'x \in\p,\] and
$dx=2^{-m} |r|^{(n-m-1)/2} dr dv$.
\end{lemma}

The following statement is new and suggestive. It contains  a
matrix generalization of the relevant   formula by Smith and
Solmon \cite[Lemma 2.2]{SS} corresponding to  the case  $m=1$.
\begin{lemma}\label{lSS}
 Let $1\leq k\leq n-m$. Then \be\label{ss} \intl_{\Ma}f(x)
dx=\frac{\sig_{n,m}}{
\sig_{n-k,m}}\intl_{\vnk}d_\ast\xi\intl_{\Mt} f(\xi z)|z|^k_m
\;dz. \ee
\end{lemma}
\begin{proof}
Let $I=\int_{\Ma}f(x) dx$, $d=(m+1)/2$. By Lemma \ref{l2.3}, \bea
\nonumber I&=& 2^{-m} \intl_{\vnm}dv \intl_{\p}
f(vr^{1/2})|r|^{n/2-d} dr\\\nonumber &=&\frac{\sig_{n,m}}{2^m}
\intl_{SO(n)}d\g \intl_{\p} f(\g v_ 0 r^{1/2})|r|^{n/2-d} dr,
\qquad \forall v_0\in\vnm. \eea Choose $v_0=\xi_0 u_0$, where
$$
\xi_0= \! \left[\begin{array} {c}  0 \\
I_{n-k}
\end{array} \right] \! \in \vnk, \qquad u_0= \! \left[\begin{array} {c}  0 \\
I_{m}
\end{array} \right]  \! \in V_{n-k,m}.
$$
By setting $\xi=\g\xi_0$, we obtain \bea \nonumber
I&=&\frac{\sig_{n,m}}{2^m} \intl_{\vnk}d_\ast\xi \intl_{\p} f(\xi
u_0 r^{1/2})|r|^{n/2-d} dr\qquad (\xi\to\xi\b) \\\nonumber&=&
\frac{\sig_{n,m}}{2^m} \intl_{\vnk}d_\ast\xi \intl_{O(n-k)} d\b
\intl_{\p} f(\xi \b u_0 r^{1/2})|r|^{n/2-d} dr\\\nonumber
&=&\frac{\sig_{n,m}}{2^m \sig_{n-k,m}} \intl_{\vnk}d_\ast\xi
\intl_{V_{n-k,m}} du  \intl_{\p} f(\xi  u r^{1/2})|r|^{n/2-d} dr
\\\nonumber &=&\frac{\sig_{n,m}}{
\sig_{n-k,m}}\intl_{\vnk}d_\ast\xi\intl_{\Mt} f(\xi z)|z|^k_m \;
dz. \eea
\end{proof}

\subsection{Riesz potentials}\label{s2.2}

We recall basic facts from  \cite{OR} and \cite{Ru6}  related to
 Riesz potentials of functions of matrix argument. These potentials  arise in different aspects of
 analysis \cite{Ge}, \cite{Kh}, \cite{St1}. They have a number of
 specific higher rank features and coincide for
  $m=1$ with classical integrals of Marcel Riesz \cite{Ru1}, \cite{SKM},
  \cite{St2}.
 In the following, we assume $m \ge 2$. The Riesz potential of
order $\a \in \bbc$ of a function $f\in\S(\Ma)$ is defined as
analytic continuation of the integral \be\label{rie} (I^\a
f)(x)=\frac{1}{\gam_{n,m} (\a)} \intl_{\Ma} f(x-y) |y|^{\a-n}_m
dy, \ee where $|y|_m =\det (y'y)^{1/2}$, \be\label{gam} \gam_{n,m}
(\a)=\frac{2^{\a m} \, \pi^{nm/2}\, \Gam_m (\a/2)}{\Gam_m
((n-\a)/2)},\quad \a\neq n-m+1, \,  n-m+2, \ldots ,\ee cf. (\ref
{ri}). This integral converges absolutely if and only if
$Re\,\a>m-1$ and extends to all $\a \in \bbc$ as a meromorphic
function whose only poles are at the points $\a=n-m+1, n-m+2,
\ldots\;$. The order of these poles is the same as in $\Gam_m
((n-\a)/2)$.
\begin{theorem}
If $f$ and $\phi$ are Schwartz functions on $\Ma$, then for all
complex $\a\neq n-m+1, \,  n-m+2, \ldots\;$, \be\label{FIa}
 (I^\a
f, \phi)=(2\pi)^{-nm}(|y|_m^{-\a}(\F f)(y),\; (\F\phi)(y)),\ee
the
expression on  each side being understood in the sense of analytic
continuation.
\end{theorem}

This statement is a consequence of the relevant functional
equation for the corresponding zeta distributions, see 
\cite{Ru6}, \cite{FK}.

If $\a=k$, $k=0,1,2,\dots , m-1$ and $k \neq    n - m + 1,
n - m + 2, \ldots $, then  $I^\a f$ is
 a convolution with a
positive measure supported by the manifold of all matrices $x$ of
rank $\leq k$. Combining this fact with the case $Re\,\a>m-1$, we
introduce the Wallach-like set \be\label{wa} \V= {\bf W}_1 \cap
{\bf W}_2, \ee  where
\[ {\bf W}_1 =  \{0,
1, 2, \ldots,  k_0 \}, \qquad   k_0=\min(m-1, n-m); \]
\[{\bf W}_2  =  \{\a : Re \, \a >  m  -  1; \; \a   \neq    n - m + 1,
n - m + 2, \ldots \}. \]   An analog of (\ref{wa}) is defined in
\cite[p. 137]{FK}  for distributions of different type. For
$\a\in\V$, one can write  \be\label{def-Rp} (I^\a f)(x)=(f\ast
\mu_\a)(x)=\intl_{\Ma} f(x-y) d\mu_\a(y), \ee \cite[Theorem
3.14]{OR}, \cite[Theorem 5.1]{Ru6}, where the measure $\mu_\a$ is
defined by 
\be\label{any} (\mu_\a,\psi) = \left \{
\begin{array} {ll}\displaystyle{\frac{1}{\gam_{n,m} (\a)} \intl_{\Ma}
|y|^{\a-n}_m \overline{\psi(y)} dy \quad
 \mbox{  if   $Re \, \a > m-1 $}}, \\
{} \\
\displaystyle{c_k\intl_{\frM_{k,m}}dy \intl_{O(n)} \overline{\psi \left (\gam
\left[\begin{array} {c} y \\ 0
\end{array} \right]  \right )} \, d\gam \quad
\mbox{if   $\a=k$}}. \\
\end{array}
\right.
 \ee
 Here, $k=1,2,\dots , n-m$, $\psi$ is a compactly supported continuous function, and  \be\label{con1} c_k=2^{-km}
\,\pi^{-km/2} \, \Gam_m\Big ( \frac{n-k}{2}\Big ) / \Gam_m\Big (
\frac{n}{2}\Big ). \ee  Owing to (\ref{FIa}), for $\a=0$, $\mu_\a$
is the usual delta function, and we set $I^0 f=f$. Note that the
sets of $\a$ in both lines of (\ref{any}) may overlap.  In  this  case
 we have two different representations of
$(\mu_\a,\psi)$.

One can use (\ref{def-Rp}) as a definition of $I^\a f$, $\a\in\V$,
for arbitrary locally integrable function provided the integral
$f\ast \mu_\a$ converges absolutely.

\begin{theorem}[\cite{Ru6}, Section 5.3]\label{t-exist}
Let $f\in\lp$, $1 \le p <n/(Re \, \a +m-1)$.

\noindent \rm{(i)}
 If $Re \, \a >  m  -  1$, then  \be\label{ine} \Big | \, \intl_{\Ma}\!\!\! \exp (-\tr
(x'x)) \, (I^\a f)(x) \, dx \Big |\le c_1 \, ||f||_p \; .\ee

\noindent \rm{(ii)}
 If $\a =k$, $k=1,2,\dots, n-m$, then
\be  \label{rav1} \intl_{\Ma} \frac{|(I^k f)(x)|}{|I_m +
x'x|^{\lam/2}} \, dx \le c_2 ||f||_p\;, \quad \ee provided $$ \lam
>k+\max \left (m-1, \frac{n+m-1}{p'}\right ), \qquad \frac{1}{p}+
\frac{1}{p'}=1. $$
\end{theorem}

This statement  shows that for $\a\in\V$ and $f\in L^p$, the
definition (\ref{def-Rp}) is meaningful, i.e., $(I^\a f)(x)$ is
finite for almost all $x\in\Ma$ provided $1 \le p <n/(Re \, \a
+m-1)$. The last equality agrees with the classical one $1 \le p <n/Re
\, \a $ for $m=1$ \cite{St2} which is sharp.  We do not know
whether the restriction $p <n/(Re \, \a +m-1)$  is necessary if
$m>1$.

\subsection{Radon transforms on the space of matrices}
The main references for this subsection are \cite{OR1}, \cite{OR},
\cite{Pe}. We fix  positive integers $k,n$, and $m$, $0<k<n$, and let 
$\vnk$ be the Stiefel manifold of orthonormal $(n-k)$-frames in
$\bbr^n$. For $\; \xi\in\vnk$ and $t\in\Mt$, the linear manifold
\be\label{plane} \tau= \tau(\xi,t)=\{x\in\Ma:\eq\} \ee
 will be called a {\it  matrix $k$-plane} in $\Ma$.  We denote by  $\Gr$ the 
 set of all such
 planes.  Each $\t \in \Gr$ is
an ordinary $km$-dimensional plane in $\bbr^{nm}$, but the set
$\frT$ has measure zero in the manifold of all such planes.

 Note that $\tau(\xi,t)=\tau(\xi\theta ', \theta 
t)$ for all $\theta\in O(n-k)$.  We identify  functions $\varphi(\t)$ on
$\Gr$ with functions $\fc$ on $\cd$ satisfying $\varphi(\xi\theta
',\theta t)=\fc$ for all $\theta\in O(n-k)$, and supply  $\Gr$
with the measure $d\tau$ so that \be \intl_{\Gr} \varphi(\t) \,
d\tau=\intl_{\cd}
 \fc \, d\xi dt. \ee

The  {\it matrix $k$-plane Radon transform} $f(x)\to\hat f(\t)$
assigns to a function $f(x)$ on $\Ma$ a collection of integrals of
$f$ over all matrix planes $ \tau \in \Gr$.
 Namely, \[ \hat f (\tau)=\int_{x \in \tau} f(x).\]  Precise meaning of this integral is
 the following:
\be\label{4.9} \hat f (\tau) \equiv \rf=\intl_{\Mkm} f\left(g_\xi
\left[\begin{array} {c} \om \\t
\end{array} \right]\right)d\om,
\ee where $ g_\xi \in SO(n)$ is a rotation satisfying
\be\label{4.24}
g_\xi\xi_0=\xi, \qquad \xi_0=\left[\begin{array} {c}  0 \\
I_{n-k} \end{array} \right] \in \vnk. \ee

The corresponding dual Radon transform $\vp(\t)\to\check \vp (x)$
assigns to a function $\vp(\tau)$ on $\Gr$ its mean value over all
matrix planes $\tau$ through $x$:
$$
\df=\intl_{\tau\ni x}\vp(\tau), \qquad x\in\Ma.
$$
This means that \be\label{4.2} \qquad \df=\intl_{\vnk}
\varphi(\xi,\xi'x)d_\ast\xi. \ee The corresponding  duality
relation reads \be\label{4.3} \intl_{\Ma} f(x)\df dx=
\intl_{\vnk}d_\ast\xi\intl_{\Mt}\fc\rf dt. \ee

\begin{theorem}[\cite{OR}]\label{t5.1} \hfill

\noindent {\rm (i)} The Radon transform $\rf$,  $f\in\lp$,  is
finite for almost all $(\xi, t) \in \cd$ if and only if
\be\label{lp} 1\leq p<p_0=\frac{n+m-1}{k+m-1}. \ee

\noindent {\rm (ii)} If $\fc$ is a  locally integrable function
on the set $\cd$, $1\leq k\le n-m$, then the dual Radon transform
$\df$ is finite
 for almost all $x \in \Ma$.
\end{theorem}

 The following statement is a matrix generalization of
the so-called projection-slice theorem. It links together the
Fourier transform (\ref{ft}) and the Radon transform (\ref{4.9}).
In the case $m=1$, this theorem can be found in [Na, p. 11] (for
$k=n-1$) and [Ke, p. 283] (for any $0<k<n$).

For $y=[y_1 \dots y_m]\in\Ma$, let $\L (y)=\span (y_1, \dots, y_m)$
be the span of the $n$-vectors $y_1, \dots, y_m$. Suppose
that $\rank (y)=\ell$. Then $\dim\L (y)=\ell\leq m$.

\begin{theorem}[\cite{Sh1}, \cite{Sh2}, \cite{OR}]\label{CST}
Let $f\in L^1(\Ma), \;1\leq k\leq n-m$. If  $y \in\Ma$, and
$\zeta$ is an $(n-k)$-dimensional plane in $\bbr^n$ containing  $\L (y)$, then for
any orthonormal frame $\xi\in\vnk$ spanning
$\zeta$, there exists $b\in\Mt$ so that $y=\xi b$. In this case
\be\label{4.20} (\F f)(\xi b)=\tilde \F[\hat f(\xi,\cdot)](b),
\quad \xi\in\vnk, \quad b\in\Mt, \ee where $\tilde \F$ stands for
the Fourier transform on $\Mt$ in the $(\cdot)$-variable.
 \end{theorem}

 \begin{corollary}[\cite{OR}]  The Radon transform
$f\rightarrow\hat f$ is injective on the Schwartz space $\S(\Ma)$
if and only if  $1\le k\le n-m$.
\end{corollary}

\section{Continuous wavelet transforms}

\subsection{Some heuristics} Following the philosophy which was  described in Introduction for the
rank-one case,  we will introduce continuous  wavelet
transforms on $\Ma$ associated to the Riesz potential (\ref{rie}).
The heuristic argument presented below shows that these  ``higher
rank'' wavelet transforms are essentially multiscaled, with the
scaling parameter
 represented by a positive definite matrix, rather then    a positive number as in the rank-one case.

We recall the notation (\ref{2.1}) for the invariant measure
$d_\ast r$ on $\p$ and  start with the following simple
observation.
\begin{lemma}\label{smm}
Let $w_0$ be a symmetric function on $\p$  satisfying
\be\label{w0} \intl_{\p}\frac{|w_0(r)|}{|r|^{(\a-n)/2}}\; d_\ast
r<\infty \quad \mbox{and}\quad
c_\a=\intl_{\p}\frac{w_0(r)}{|r|^{(\a-n)/2}}\;d_\ast r\not =0. \ee
 Then for  $s\in\p$, \be\label{xm-w}
|s|^{(\a-n)/2}=c_\a^{-1}
\intl_{\p}\frac{w_0(a^{-1/2}sa^{-1/2})}{|a|^{(n-\a)/2}}\;d_\ast a
\; .\ee
\end{lemma}

\begin{proof} Using the symmetry (\ref{sym}) and
  changing variable $a=\rho^{-1}$,
$d_\ast a=d_\ast\rho$, we rewrite  (\ref{xm-w}) as
$$
|s|^{(\a-n)/2}=c_\a^{-1}\intl_{\p}\frac{w_0(s^{1/2}\rho
s^{1/2})}{|\rho|^{(\a-n)/2}}\;d_\ast\rho\, .
$$
It remains to set $s^{1/2}\rho s^{1/2}=r$.
\end{proof}

According to (\ref{xm-w}), for $Re\,\a
>m-1$, the Riesz potential (\ref{rie})  is
represented as \bea\nonumber  (I^\a
f)(x)\!\!\!&=&\!\!\!\frac{c_\a^{-1}}{\g_{n,m}(\a)}\intl_{\p}|a|^{(\a-n)/2}
d_\ast a  \intl_{\Ma} \!\!\!f(x-y) w_0(a^{-1/2} y'y
a^{-1/2})dy\\\label{eq8}
&=&\frac{c_\a^{-1}}{\g_{n,m}(\a)}\intl_{\p}|a|^{\a/2} d_\ast a
\intl_{\Ma} f(x-ya^{1/2}) w_0( y'y )dy. \eea The inner integrals
in these expressions resemble the  wavelet transform (\ref{wav})
and inspire the following.
\begin{definition}\label{def-wt}
 Let $w(y)=w_0(y'y)$, $y\in\Ma$, be a radial
function satisfying certain cancellation conditions (which depend on the context). Let
$w_a(y)=|a|^{-n/2}w(ya^{-1/2})$, $a\in\p$.  We call
\bea\label{defW} (\W_a f)(x)&=&\intl_{\Ma} f(x-ya^{1/2})\,  w(y) dy \\
&=&(f\ast w_a)(x),\qquad x\in\Ma, \nonumber \eea the {\it
continuous wavelet transform of} $f$
 generated by the wavelet function $w$ and the $\p$-valued scaling parameter $a$.
\end{definition}

Note that the symmetry condition for  $w_0$ indicated  in Lemma
\ref{smm} plays an auxiliary  role. It was imposed only for
technical reasons and not included in Definition \ref{def-wt}. In
the following, this condition will appear  on the Fourier
transform side for $\F w$.

If the function $w_0$ in Definition \ref{def-wt}
  is  symmetric,
then one can write (\ref{eq8}) as \be\label{Ia-W} (I^\a f)(x) =
c_{n,m}(\a,w) \intl_{\p} (\W_a f)(x) |a|^{\a/2}\, d_\ast a, \;
\ee
$$
Re \, \a > m-1; \quad \a \neq n-m+1, n-m+2, \ldots,
$$
$c_{n,m}(\a,w)\equiv \const$. This formula is expected to be true
for other values of $\a$  (e.g., for $\a=0$) if $w$ obeys certain
cancellation conditions. Of course, if $ Re \, \a \le m-1$, then 
the integral on the right-hand
side of (\ref{Ia-W})  diverges in general, and must be interpreted in a
suitable way depending on a class of functions $f$ and a choice of  the wavelet $w$.

One can replace the function  $w$ in (\ref{defW})  by a
finite  radial  Borel measure $\nu$  on $\Ma$ so that
$$\intl_{\Ma} \psi (\g x) \, d\nu(x)=\intl_{\Ma} \psi (x)\,  d\nu(x)$$
for all $\g\in SO(n)$ and $ \psi\in C_c(\Ma)$. If $\nu$ obeys 
some cancellation (for instance, $(\F \nu)(y)\equiv 0$ on matrices
$y$  of rank $<m$) we call \be\label{defWm} (\W_{\nu,\; a}
f)(x)\equiv(f\ast \nu_a)(x) = \intl_{\Ma} f(x-ya^{1/2}) d\nu(y)
\ee  the  wavelet transform of $f$ generated by the wavelet
measure $\nu$.

\subsection{Calder\'{o}n's reproducing formula}
Denote formally \be I(\nu, f)(x)= \intl_{\p}(\W_{\nu,\; a} f)(x)\,
d_\ast a, \ee where, as above, $d_\ast a= |a|^{-d} da$,
$d=(m+1)/2$.
 The following statement justifies (\ref{Ia-W}) for $\a=0$ and
generalizes the classical Calder\'{o}n reproducing formula (cf.
Theorem 1 in \cite{Ru3} ) to  functions of matrix argument.

\begin{theorem}\label{tCald}
Let  $\nu$ be a radial finite Borel measure on $\Ma$ so that for
all $y \in \Ma$ of rank $m$,  we have $(\F \nu)(y)=u_0(y'y)$, where
$u_0(r)$ is a symmetric function on $\p$. Suppose that the
integral \be c_\nu= \lim\limits_{A\to 0 \atop  B\to \infty} \,
\frac{2^m}{\sig_{n,m}} \intl_{\{y\in\Ma:\,A<y'y<B\}}\frac{(\F
\nu)(z)}{|z|^{n}_m}\, dz\quad \ee ($A,B \in \p$) is finite. Then
for $f\in L^2(\Ma)$, \be c_\nu  f(x)=I(\nu, f)(x)\equiv
\underset{\rho\to \infty} {\lim\limits_{\e\to 0}^{(L^2)}} \,
\intl_{\e I_m}^{\rho I_m}(\W_{\nu,\; a} f)(x)\, d_\ast a. \ee
\end{theorem}

\begin{proof}
For $0<\e<\rho<\infty$, let  \be I_{\e,\rho}(\nu,f)(x)= \intl_{\e
I_m}^{\rho I_m}(\W_{\nu,\; a} f)(x)\, d_\ast a \ee and assume
first that $f\in L^1\cap L^2$. Then, by the generalized Minkowski
inequality, $I_{\e,\rho}(\nu,f)\in L^1\cap L^2$, and we have
\be\label{eq7} \F [I_{\e,\rho}(\nu,f)](y)= \intl_{\e I_m}^{\rho
I_m}\F[\W_{\nu,\; a} f](y)\, d_\ast a.\ee By taking into account
that $ \F[\W_{\nu,a} f](y)=(\F f)(y)(\F \nu_a)(y),$ and \be
\nonumber (\F \nu_a)(y)= \intl_{\Ma} \exp(\tr(ia^{1/2}y'x))
d\nu(x)=(\F \nu)(ya^{1/2}),\ee we obtain \be \label{FIer}\F
[I_{\e,\rho}(\nu,f)](y)=k_{\e,\rho}(y)(\F f)(y),\ee where
$$ k_{\e,\rho}(y)= \intl_{\e I_m}^{\rho I_m}(\F
\nu)(ya^{1/2})\, d_\ast a =\intl_{\e I_m}^{\rho
I_m}u_0(a^{1/2}y'ya^{1/2})\, d_\ast a.$$ Suppose that $\rank(y)=m$
(the set of all such $y$ has a full measure in $\Ma$). Since $u_0$ is
symmetric, then $u_0(a^{1/2}ra^{1/2})=u_0(r^{1/2}a r^{1/2})$,
$r=y'y$, and the change of variable $s=r^{1/2}ar^{1/2}$,
 yields \bea\nonumber k_{\e,\rho}(y)&=&\intl_{\e
I_m}^{\rho I_m }u_0(r^{1/2}ar^{1/2})\,d_\ast a= \intl_{\e r}^{\rho
r } u_0(s)\,d_\ast s\qquad \mbox{(use Lemma \ref{l2.3}) }
\\&=&\frac{2^m}{\sig_{n,m}}\intl_{\{z\in\Ma:\,\e r<z'z<\rho r\}}\frac{(\F
\nu)(z)}{|z|^{n}_m}\, dz,\qquad r=y'y.\nonumber\eea Since
$k_{\e,\rho}(y)$ is bounded uniformly in $\e, \rho$, and $y$, then
by the Lebesgue theorem on dominated convergence, \be\label{eq1}
\|I_{\e,\rho}(\nu, f)-c_\nu
f\|_2=(2\pi)^{-nm/2}\|(k_{\e,\rho}-c_\nu )\F f\|_2\to 0 \ee as
$\e\to 0$, $\rho\to\infty$. This proves the statement for $f\in
L^1\cap L^2$. A standard procedure allows us to extend the result
to all $f\in L^2$. We recall this argument for convenience of the
reader. For any $f\in L^2$, we have  $$ \|I_{\e,\rho}(\nu,
f)\|_2\leq \intl_{\e I_m}^{\rho I_m}\|\W_{\nu,\; a} f\|_2 \,
d_\ast a\leq c_{\e,\rho}\|f\|_2\|\nu\|=c_{\e,\rho,\nu} \|f\|_2
$$
where $||\nu ||$ stands for the total variation of  $|\nu|$, $
c_{\e,\rho}=\const$, $c_{\e,\rho,\nu}=c_{\e,\rho}\,||\nu ||$.
Given a small $\del >0$,  we choose   $g\in L^1\cap L^2$ so that
$\|f-g\|_2<\del$. Since $k_{\e,\rho}$ is uniformly bounded, then
(\ref{FIer}) (with $f$ replaced by $g$) implies the uniform
estimate
$$
 \|I_{\e,\rho}(\nu, g)\|_2 \leq c \|g\|_2 \leq c \|g-f\|_2+c
 \|f\|_2,
$$
and, therefore,
 \bea\nonumber \|I_{\e,\rho}(\nu, f)\|_2&\leq&
\|I_{\e,\rho}(\nu, f-g)\|_2+\|I_{\e,\rho}(\nu,
g)\|_2\\\nonumber&\leq& \del c_{\e,\rho,\nu}+ c \|f-g\|_2+c
\|f\|_2\\\nonumber&\leq& \del (c_{\e,\rho,\nu}+ c) +c \|f\|_2.
\eea Assuming  $\del\to 0$, we obtain  $\|I_{\e,\rho}(\nu,
f)\|_2\leq c \|f\|_2$. This gives \bea\nonumber \|I_{\e,\rho}(\nu,
f)-c_\nu f\|_2&\leq& \|I_{\e,\rho}(\nu,
f-g)\|_2+\|I_{\e,\rho}(\nu,g)-c_\nu g\|_2
\\\nonumber&+&c_\nu\|g-f\|_2 \\\nonumber &\leq&(c+ c_\nu)\del
+\|I_{\e,\rho}(\nu,g)-c_\nu g\|_2.\eea Owing to (\ref{eq1}) (with
$f$ replaced by $g$),
$$
\|I_{\e,\rho}(\nu, f)-c_\nu f\|_2\to (c+ c_\nu)\del \quad
\mbox{as} \quad \e\to 0,\; \rho\to\infty.$$ Since $\del$ is
arbitrarily  small, we are done.
\end{proof}

\section{Inversion of  Riesz potentials}\label{s4}

\setcounter{equation}{0}
\noindent
 We recall that the wavelet transform of a function $f$ on $\Ma$
 is defined by
 $$(\W_a
f)(x)=\intl_{\Ma} f(x-ya^{1/2})\,  w(y)\,  dy=(f\ast w_a)(x)$$
where \be \label{w-fun}w_a(x)=|a|^{-n/2}w(xa^{-1/2}), \qquad
x\in\Ma, \quad a\in\p. \ee Owing to (\ref{FIa}), it is natural to
 expect, that the inverse of the Riesz potential (\ref{def-Rp})
can be  obtained if we formally replace $\a$ by $-\a$ in
(\ref{Ia-W}); cf. (\ref{inv}).  This gives \be\label{Da}
f(x)=c_{n,m}(-\a,w)\intl_{\p} \frac{(\W_a I^\a f)(x)}{
|a|^{\a/2}}\, d_\ast a. \ee Below we give this
  formula precise meaning.
\begin{theorem}\label{in2}
Let  $\a\in\V$ and $ f\in L^2\cap L^p $ for some $p$ satisfying
$$  1 \le p <\frac{n}{Re\,\a +m-1}\; .$$
Suppose that  $w \in \S(\Ma)$ is a   radial function  such that

\noindent
{\rm (a)}  $(\F w)(y)=u_0(y'y)$, where $u_0(r)$ is a
symmetric function on $\p$ vanishing identically  in a neighborhood of the
boundary $\partial\p$;

\noindent {\rm (b)}  The integral \bea\label{dwa}
d_w(\a)&=&\frac{2^m}{\sig_{n,m}} \intl_{\Ma}\frac{(\F
w)(z)}{|z|^{n+\a}_m}\, dz\\\nonumber&=&\lim\limits_{ B\to
\infty}\frac{2^m}{\sig_{n,m}} \intl_{\{z\in\Ma:\,
z'z<B\}}\frac{(\F w)(z)}{|z|^{n+\a}_m}\, dz \eea $(B \in \p)$ is
finite. Then \bea \label{invR}d_w (\a) f&=&\intl_{\p} \frac{\W_a
I^\a f}{ |a|^{\a/2}}\, d_\ast a = \lim\limits_{\e \to 0\atop
\rho\to\infty }^{(L^2)} \intl_{\e I_m}^{\rho I_m} \frac{\W_a I^\a
f}{ |a|^{\a/2}}\, d_\ast a. \eea
\end{theorem}

\begin{proof}
We observe that the integrand in (\ref{dwa}) has no singularity at
$|z|=0$ thanks to the assumption (a) above. Moreover,  for $Re \,
\a
>m-1$, this integral  is  finite automatically,  because
$$
\intl_{\Ma}\frac{(\F w)(y)}{|y|^{n+\a}_m}\, dy \le c
\intl_{\Ma}|I_m +y'y|^{-(n+\a)/2} dy< \infty, \quad c \equiv
c(w),$$
 see formula (A.6) in \cite{OR}.

To prove the theorem, we set \be\label{Dae} (T^\a_{\e,\rho}
\vp)(x)=\intl_{\e I_m}^{\rho I_m } \frac{(\W_a \vp)(x)}{
|a|^{\a/2}}\, d_\ast a,\qquad 0<\e<\rho<\infty  \ee (if  $Re \, \a
>m-1$ one can assume $\rho=\infty$), and show that \be
\label{FDI} T^\a_{\e,\rho} I^\a f=\F^{-1}[\psi_{\e,\rho}^\a \F f]
, \ee \be \psi_{\e,\rho}^\a(y)=\frac{2^m}{\sig_{n,m}}
\intl_{\{z\in\Ma:\,\e(y'y)<z'z<\rho(y'y)\}}\frac{(\F
w)(z)}{|z|^{n+\a}_m}\, dz. \ee
As we have shown this,  by the Lebesgue  dominated convergence
theorem, owing to (\ref{dwa}) and the
uniform boundedness of $\psi_{\e,\rho}^\a$, we obtain the desired
result:
$$
\|T^\a_{\e,\rho} I^\a f-
d_w(\a)f\|_2=(2\pi)^{-mn/2}\|(\psi_{\e,\rho}^\a - d_w(\a))\F f\|_2
\to 0 \quad $$ as  $\e\to 0,\; \rho\to\infty.$

We observe that  for any $f\in L^p$ and $w\in L^1$,  \be \label
{DeI}T^\a_{\e,\rho} I^\a f=I^\a g , \qquad g=T^\a_{\e,\rho} f.\ee
The validity of interchange of integrals follows  by Theorem
\ref{t-exist}, according to which, the integral $I^{|\a|} \,
[|f|\ast |w_a| ](x)$ is finite for almost all $x$ because $
|f|\ast |w_a|\in L^p $.

 We first prove (\ref{FDI}) for $f$
belonging to the Schwartz space $\S=\S(\Ma)$. By (\ref{defW}),
$$ \F(W_a f)(y)=\F(f\ast w_a)(y)=(\F f)(y)(\F w_a)(y),$$ where
\bea \nonumber (\F
w_a)(y)&=&|a|^{-n/2}\intl_{\Ma} \exp(\tr(iy'x)) w (xa^{-1/2}) dx\\
\nonumber&=& \intl_{\Ma} \exp(\tr(ia^{1/2}y'z)) w (z) dz=(\F
w)(ya^{1/2}).\eea Hence, $$ (\F g)(y)= h_{\e,\rho}(y) (\F
f)(y),\qquad h_{\e,\rho}(y)=\intl_{\e I_m}^{\rho I_m } \frac{(\F
w)(ya^{1/2})}{ |a|^{\a/2}}\, d_\ast a. $$ By taking into account
that  $(\F w)(y)=u_0(y'y)$ where $u_0(r)$ is symmetric, we obtain
\bea h_{\e,\rho}(y)&=&\intl_{\e I_m}^{\rho I_m }
\frac{u_0(a^{1/2}y'ya^{1/2})}{ |a|^{\a/2}}\,d_\ast a \nonumber \\
&=& \frac{2^m|y'y|^{\a/2}}{\sig_{n,m}}\intl_{\{z\in\Ma:\,\e
(y'y)<z'z<\rho(y'y)\}}\frac{(\F w)(z)}{|z|^{n+\a}_m}\, dz\nonumber
\\&=&|y|^\a_m \psi_{\e,\rho}^\a (y)\eea (see the argument  in the
proof of Theorem \ref{tCald}). This gives \be\label{Fg}(\F
g)(y)=|y|_m^{\a}\psi_{\e,\rho}^\a(y)(\F f)(y). \ee Since  $w\in\S$
and  $u_0$ is supported away  from the boundary $\partial\p$, it
follows that  $(\F g)(y) \equiv h_{\e,\rho}(y)(\F w)(y)\in\S$, and
therefore,   $g\in\S$. Hence, by (\ref{DeI}), (\ref{FIa}), and
(\ref{Fg}), for any compactly supported $C^\infty$ function $\phi$, 
we have \bea\nonumber
(T^\a_{\e,\rho} I^\a f,\;\phi)&=&(I^\a g,\; \phi)\\
\nonumber &=&(2\pi)^{-nm}(\F I^\a g,\, \F
\phi)\\
\nonumber &=& (2\pi)^{-nm}(|y|_m^{-\a}(\F g)(y),\, (\F
\phi)(y))\\\nonumber &=&(2\pi)^{-nm}( \psi_{\e,\rho}^\a (y)(\F
f)(y),\, (\F \phi)(y)).\eea
Thus, by the Parseval equality,
 \be\label{DeIa} (T^\a_{\e,\rho} I^\a
f,\;\phi)=(\F^{-1}[ \psi_{\e,\rho}^\a  \F f],\, \phi)\ee (note
that $\psi_{\e,\rho}^\a (y)(\F f)(y)\in L^2$ in view of the
boundedness of $\psi_{\e,\rho}^\a (y)$). Since
$\F^{-1}[\psi_{\e,\rho}^\a \F f]\in L^2$ and $T^\a_{\e,\rho} I^\a
f=I^\a g$ is a locally integrable function (see (\ref{ine}) and
(\ref{rav1})), then (\ref{DeIa}) implies the pointwise equality
(\ref{FDI}) for any $f\in \S$.

To complete the proof, it remains to  extend (\ref{FDI}) to all
$f\in L^2\cap L^p$. Following Theorem \ref{t-exist}, we introduce
the weighted space
$$
X=\{\vp: \;\|\vp\|_X=\intl_{\Ma} |\vp(x)| \, \om (x)\, dx <\infty
\}
$$ where $\om (x)= \exp (-\tr (x'x))$ if  $Re \, \a >m-1$, and
 $\om (x)=  |I_m +
x'x|^{-\lam/2}$ if $\a=k, \; k=1,2, \ldots , n-m$; see (\ref{rav1}). It may happened that these domains of $\a$
 overlap, but this is not important. By H\"{o}lder's inequality,
$||\vp ||_X  \le c \,||\vp ||_2 $. Since by  Theorem \ref{t-exist},
$$
\|T^\a_{\e,\rho} I^\a
f\|_X=\| I^\a g\|_X \leq c \,  \| g \|_p = c \,  \|T^\a_{\e,\rho} f\|_p
 \le c_{\e,\rho}\nf,
$$
and
$$
\|\F^{-1}[ \psi_{\e,\rho}^\a  \F f]\|_X\leq c \|\F^{-1}[
\psi_{\e,\rho}^\a  \F f]\|_2 \leq
 c'_{\e,\rho}\|f\|_2,
$$
operators $$ T^\a_{\e,\rho} I^\a : L^p \to X, \qquad \F^{-1}
\psi_{\e,\rho}^\a  \F :
 L^2 \to X
$$
are bounded. This remark allows us to  extend (\ref{FDI}) to all
$f\in L^2\cap L^p$ by taking into account that there is a sequence
$\{f_j\} \subset \S$ such that the quantities $\|f-f_j\|_p $ and
$\|f-f_j\|_2$ tend to $0$ as $j \to\infty$ simultaneously.  Such a
sequence can be explicitly constructed  using the standard
 ``averaging-truncating'' procedure.
\end{proof}

\section{Continuous ridgelet
transforms and  inversion of the Radon transform}
\setcounter{equation}{0}

\subsection{Intertwining operators }\label{ss4} Given a sufficiently good
function $w$ on $\Mt$, consider the intertwining operator
\bea\label{wxt} (W f)(\xi, t)&=& \intl_{\Ma} f(x)\,w(t-\xi '
x)\,dx \eea which transforms a function $f$ on $\Ma$ into a
function $Wf$ on the ``cylinder" $\vnk \times \Mt$. The
corresponding dual operator is defined by \bea\label{wdxt}
(W^{\ast}\vp)(x)&=&\intl_{\vnk}\!\!d_\ast
\xi\!\!\intl_{\Mt}\!\!\!\!\!\!\vp(\xi,t)\,w( \xi ' x-t)\, dt,\eea
so that \be\label{dualw} \intl_{\Ma} f(x)(W^\ast\vp)(x) dx=
\intl_{\vnk}d_\ast\xi\intl_{\Mt}\fc (Wf)(\xi,t) dt \ee (at least,
formally).  We shall see that $\p$-scaled versions of $W$ and
$W^\ast$ can be regarded as matrix modifications of continuous
$k$-plane  ridgelet transforms (see \cite{Ca}, \cite{Ru5},  and references
therein), and used for explicit and approximate inversion of the
Radon transform (\ref{4.9}). We start with some preparations.

\begin{lemma}
Given a function $\fc$ on $\cd$, let \be\label{con} (W_0\,
\vp)(\xi, t)=\intl_{\Mt} \vp(\xi,z)\,w(t- z)\, dz \ee be a
convolution in the $t$-variable. Then \be\label{w-r} (W f)(\xi,
t)=(W_0 \hat f)(\xi, t),\ee
 \be\label{dw-r} (W^{\ast}\vp)(x)= (W_0 \vp)^\vee (x)\ee
provided that either side of the corresponding equality is finite
for  $f$, $\vp$, and $w$
 replaced by $|f|$, $|\vp|$, and $|w|$, respectively.
\end{lemma}
\begin{proof}
The equality (\ref{dw-r}) follows immediately from (\ref{wdxt}).
 To prove (\ref{w-r}),
 we choose a rotation $ g_\xi \in SO(n)$  satisfying $$
g_\xi\xi_0=\xi, \qquad \xi_0=\left[\begin{array} {c}  0 \\
I_{n-k} \end{array} \right] \in \vnk. $$ The change of variable
$x=g_\xi y$ in (\ref{wxt}) gives
$$
(W f)(\xi, t)= \intl_{\Ma} f(g_\xi y)\,w( t-\xi_0 ' y)\,dy.
$$
By setting
 \[ y=\left[\begin{array} {c} \om \\ z
\end{array} \right], \qquad \om \in \frM_{k,m}, \qquad z \in
\frM_{n-k,m},\] so that $\xi_0 ' y=z$, owing to (\ref{4.9}),  we
obtain \bea \nonumber(W f)(\xi, t)&=& \intl_{\frM_{n-k,m}}w(
t-z)\; dz\intl_{\frM_{k,m}} f\left(g_\xi \left[\begin{array} {c}
\om \\ z
\end{array} \right]\right)d\om\\\nonumber
&=&\intl_{\frM_{n-k,m}}\hat f(\xi, z)\,w( t-z)\; dz=(W_0 \hat
f)(\xi, t). \eea
\end{proof}

\begin{lemma} \label{loc}Let $f(x)$ and $w(z)$ be  integrable  functions  on $\Ma$ and
$\Mt$, respectively. Then $(W^\ast \hat f)(x)$ is a
locally integrable function on $\Ma$ which belongs to $\S'(\Ma)$ and satisfies \be\label{eq4} 
\intl_{\Ma} \! 
(W^\ast \hat f)(x) \, \phi (x) \, dx = \! \intl_{\Ma} \!  f(x) \, 
(W^\ast \hat\phi) (x) \, dx, \quad 
\phi\in\S(\Ma).\ee
\end{lemma}
\begin{proof} 
We have 
\bea \mbox{\rm l.h.s}&\stackrel{\rm
(\ref{dualw})}{=}&
\intl_{\vnk}d_\ast\xi\intl_{\Mt}\rf (W\phi)(\xi,t) dt \nonumber \\
\label{fff}&\stackrel{\rm (\ref{w-r})}{=}&
\intl_{\vnk}d_\ast\xi\intl_{\Mt}\rf (W_0 \hat \phi)(\xi, t)
dt\\ &\stackrel{\rm (\ref{4.3})}{=}&
\intl_{\Ma} \!  f(x) \,  (W_0 \hat \phi)^\vee(x)
 \, dx \nonumber \\
&\stackrel{\rm (\ref{dw-r})}{=}& \mbox{\rm r.h.s} .
\nonumber \eea 
These calculations are well justified and all
statements of the lemma become clear, owing to the following
estimate of the expression (\ref{fff}):
 \bea\label{es-wd}
&\displaystyle \intl_{\vnk}&d_\ast \xi\intl_{\Mt} |\hat
f(\xi,z)|\; dz\!\!\intl_{\Mt}\!\!\! |\hat \phi(\xi,t)\,w( t-z)|\,
dt\\\nonumber&\leq&\|\hat\phi\|_\infty \|w\|_1 \intl_{\vnk}d_\ast
\xi\intl_{\Mt} |\hat f(\xi,z)|\; dz\\\nonumber & \stackrel{\rm
(\ref{4.3})}{=}&\|\hat\phi\|_\infty \|w\|_1 \|f\|_1.
 \eea
\end{proof}

In the sequel, it is convenient to use different notations for the
Fourier transform  on $\Ma$ and  $\Mt$. For the first one we write
$\F$ as before, and the second will be denoted by $\tilde \F$.
\begin{lemma}
If $w\in L^1(\Mt)$ and $\phi\in\S(\Ma)$, then
\bea\label{eq5}\qquad  (W^\ast \hat\phi)(x)&=&(2\pi)^{(k-n)m}
\intl_{\vnk} d_\ast\xi
\\ &\times& \intl_{\Mt} \exp(-\tr(ix'\xi z))(\tilde\F
w)(z)(\F\phi)(\xi z) dz. \nonumber\eea
\end{lemma}
\begin{proof}
Since $w\in L^1(\Mt)$ and $\phi\in\S(\Ma)$, the convolution
$$
(W_0 \;\hat\phi)(\xi,t)=\intl_{\Mt}\hat\phi(\xi,t-z)\,w( z)\, dz
$$
has the Fourier transform (in the $t$-variable) belonging to
$L^1(\Mt)$. Hence (see, e.g.,  \cite[p. 11]{SW}) one can write
\bea \nonumber(W_0
\;\hat\phi)(\xi,t)&=&(2\pi)^{(k-n)m}\intl_{\Mt}\!\!\!\!\!\!
\exp(-\tr(it' z))[\tilde\F (W_0 \;\hat\phi)(\xi,\cdot)](z)
dz\\\nonumber &=&(2\pi)^{(k-n)m}\intl_{\Mt}\!\!\!\!\!\!
\exp(-\tr(it' z))(\tilde\F w)(z)[\tilde\F \hat\phi(\xi,\cdot)](z)
dz. \eea By the projection-slice theorem (see (\ref{4.20})),
$$
(W_0 \;\hat\phi)(\xi,t)=(2\pi)^{(k-n)m}\intl_{\Mt}\!\!\!\!\!\!
\exp(-\tr(it' z))(\tilde\F w)(z)(\F \phi)(\xi z) dz.
$$
This proves the statement.
\end{proof}

\subsection{Continuous ridgelet transforms} Let $w(z)$  be a
sufficiently good
 function on $\Mt$, $1\leq k\leq n-m$. We consider the $\p$-
scaled version of $w$ defined by $w_a(z)=|a|^{(k-n)/2}w(za^{-1/2})$,
 $a\in\p$,
and introduce the following dual pair of intertwining operators
\bea\label{rt}(W_a f)(\xi, t)&=& \intl_{\Ma} f(x)\,w_a(t-\xi '
x)\,dx,\\\label{rdt} (W_a^{\ast}\vp)(x)&=&\intl_{\vnk}\!\!d_\ast
\xi\!\!\intl_{\Mt}\!\!\!\!\!\!\vp(\xi,t)\,w_a(\xi ' x-t)\, dt.\eea

If $w(z)$ oscillates in a certain sense we call (\ref{rt}) the
{\it continuous ridgelet transform} of $f$, and (\ref{rdt}) the
{\it dual continuous ridgelet transform} of $\vp$.

Operators (\ref{rt}) and (\ref{rdt}) generalize  usual $k$-plane
ridgelet transforms \cite{Ca}, \cite{Ru5} to the higher rank case $m> 1$.
The function $x\to w(\xi'x-t)$ is constant on each matrix plane
$\tau=\{x\in\Ma:\eq\}$ and represents a ``plane wave''.

The following definition will be useful in the sequel.

\begin{definition}\label{adm}
A function $w(z)$ on $\Mt$ is called an {\it admissible wavelet
function} if it obeys the following conditions:

\noindent (i) $w(z)$ is radial, i.e., $w(z)\equiv w_0(z'z)$, and
belongs to $L^1(\Mt)$.

\noindent (ii) The Fourier transform of $w$ has the form
$(\tilde\F w)(y)=u_0(y'y)$, where $u_0$ satisfies the symmetry
condition (\ref{sym}).

\noindent (iii)
 The integral  \bea\label{cw}
c_w&=&\frac{2^{m(k+1)} \pi^{km}}{\sig_{n,m}}\intl_{\Mt}\frac{(\tilde\F
w)(\z)}{|\z|^{n}_m}\, d\z\\\nonumber&=&
 \lim\limits_{A\to 0\atop B\to \infty} \frac{2^{m(k+1)} \pi^{km}}{\sig_{n,m}}
\intl_{\{\z\in\Mt:\,A<\z'\z<B\}} \, \frac{(\tilde\F
w)(\z)}{|\z|^{n}_m}\, d\z \eea ($A,B\in\p$) is finite.
\end{definition}

\subsection{Inversion of the Radon transform}
\subsubsection {Discussion of the problem} There exist different
approaches to inversion of the Radon transform (\ref{4.9}); see
 \cite{OR}. The consideration below sheds new light on
this problem and provides essential progress. To explain our
strategy, we use intertwining fractional
integrals $P^{\a}f$    and $\pd\vp $ of the Semyanistyi type, which link together the Radon
transform $f(x)\to\rf$, the dual Radon transform
$\vp(\xi,t)\to\df$, and Riesz potentials. Namely, we define
 \be\label{ppd0} P^{\a}f=\tilde I^\a\hat f, \qquad
\pd\vp=(\tilde I^\a\vp)^\vee, \ee  $$\a\in\bbc,\qquad \a  \neq n-k
- m + 1, \, n-k  - m + 2,\dots .$$ Here, $1\le k \le n-m$ and
$\tilde I^\a$ denotes the Riesz potential on $\Mt$  acting in the
$t$-variable. Operators (\ref{ppd0})  were introduced in  \cite{OR}. If
$Re\,\a
>m-1$, they are represented as absolutely convergent integrals
  \bea\label{P0}
 \qquad \qquad (P^\a f)(\xi, t)\label{p0}&=&
\frac{1}{\g_{n-k,m}(\a)}\intl_{\Ma} f(x)\,|\xi ' x-t|_m^{\a+k-n}\,dx,\\
(\pd\vp)(x)\!\!&=&\!\!\frac{1}{
\g_{n-k,m}(\a)}\intl_{\vnk}\!\!d_\ast
\xi\!\!\intl_{\Mt}\!\!\!\!\!\!\vp(\xi,t)\,|\xi' x-t|_m^{\a+k-n}\,
dt\label{pd0}, \eea where $\g_{n-k,m}(\a)$ is the normalizing
constant in the definition of the Riesz potential on $\Mt$; cf.
(\ref{sem}), (\ref{gam}).

The  following statement determines our way of thinking.

\begin{theorem}[\cite{OR}, Section 5.3]\label{t0.9}
Let $1\le k\le n-m$, $ \; \a\in\V$. Suppose that
 $$
f\in L^p(\Ma), \qquad 1\leq p<\frac{n}{Re\,\a+k+m-1}\;.$$
 Then \be \label{0.5}(\pd\hat f)(x)=
c_{n,k,m} (I^{\a+k} f)(x),\ee
 \be\label{cnkm}
 c_{n,k,m}=2^{km}\pi^{km/2}\gm\left(\frac{n}{2}\right)/\gm\left(\frac{n-k}{2}\right).
 \ee
 In particular, for $\a=0$,
 \be\label{fu0}(\hat f)^{\vee} (x) \! = \! c_{n,k,m} (I^k f)(x) \ee (the
generalized Fuglede formula).
\end{theorem}

Formula (\ref{0.5}) paves two ways to the inversion of the Radon
transform. Following the first one, we use  (\ref{0.5}) as it is
and  invert the Riesz potential $I^{\a+k} f$ by choosing $\a$ as
we wish. For instance, one can set $\a=0$ and apply (\ref{fu0}).
This program  can be realized using results of the previous
section. The second way is to  set (formally)   $\a=-k$ in
(\ref{0.5}). This gives \be\label{gfu2} c_{n,k,m}
f(x)=(\stackrel{*}{P}\!{}^{-k}\hat f)(x)=(\tilde I^{-k} \hat
f)^\vee (x),\ee and we have to find ``good" representation for the
inverse of the Riesz potential $\tilde I^k$ applied to $\hat f
(\xi, t)$ in the $t$-variable. In the first case, we just
apply the left inverse operator to $I^{\a+k}$. In the second one,
we do not know in advance whether $\hat f (\xi, \cdot)$ lies in the
range of the Riesz potential $\tilde I^k$. To circumvent this
difficulty, we make use of continuous ridgelet transforms.

Below we consider both approaches.

\subsubsection {The first method}  We utilize the generalized Fuglede
formula (\ref{fu0}) and invert  the Riesz potential according to
Theorem \ref{in2}. This gives the following result for
 the Radon transform.
\begin{theorem}
Let $1\le k\leq n-m$, \be f\in L^2\cap L^p, \qquad 1\leq
p<\frac{n}{k+m-1}.\ee Suppose that $\W_a$ is the continuous
wavelet transform (\ref{defW}) generated by the wavelet $w$
satisfying conditions of Theorem \ref{in2}.  Then the Radon
transform $f(x) \to\rf$ can be inverted by the formula \be
\label{invRr1}
 d_w (k) \, f=c_{n,k,m}\intl_{\p} \frac{\W_a \hat f}{
|a|^{k/2}}\, d_\ast a= c_{n,k,m} \lim\limits_{\e \to 0\atop
\rho\to\infty }^{(L^2)} \intl_{\e I_m}^{\rho I_m} \frac{W_a \hat
f}{ |a|^{k/2}}\, d_\ast a \ee where
$$
c_{n,k,m}=\frac{2^{km}\pi^{km/2}\gm (n/2)}{\gm((n-k)/2)}, \qquad
d_w(k)=\frac{2^m}{\sig_{n,m}} \intl_{\Ma}\frac{(\F
w)(y)}{|y|^{n+k}_m}\, dy. $$
\end{theorem}

\subsubsection {The second method}
 By (\ref{invR}) and (\ref{gfu2}),
it is natural to expect,  that the Radon transform can be
inverted as
$$
f(x)=\intl_{\vnk}\!\!\!\!d_\ast\xi \intl_{\p } \frac{(\hat f (\xi,
\cdot) \ast w_a)(\xi' x)}{ |a|^{k/2}}\, d_\ast a =  \intl_{\p }
\frac{(W_a^\ast \hat f)(x)}{ |a|^{k/2}}\, d_\ast a
$$
(up to a constant multiple) where $W_a^\ast $ is the dual ridgelet
 transform (\ref{rdt}). Below we justify  this formula.

\begin{theorem}\label{tinv1}
Let $f\in L^1(\Ma)\cap L^2(\Ma)$, $1\leq k\leq n-m$. If $w$ is an
admissible wavelet function (see Definition \ref {adm}), then the
Radon transform $f(x) \to\rf$  can be inverted by the formula
\be\label{inv1} c_w \,f= \intl_{\p } \frac{W_a^\ast \hat f}{
|a|^{k/2}}\, d_\ast a = \underset {\rho\to\infty }{
\lim\limits_{\e \to 0}^{(L^2)}} \, \intl_{\e I_m}^{\rho I_m }
\frac{W_a^\ast \hat f}{ |a|^{k/2}}\, d_\ast a, \ee where $c_w$ is
defined by (\ref{cw}),  and  \be\label{Wa-r} (W_a^{\ast} \hat
f)(x)=\intl_{\vnk}\!\!d_\ast \xi \intl_{\Mt} \hat
f(\xi,t)\,w_a(\xi ' x-t)\, dt.\ee
\end{theorem}
\begin{proof}
Consider the truncated integral \be\label{Ier} I_{\e,\rho} f=
\intl_{\e I_m}^{\rho I_m } \frac{W_a^\ast \hat f}{ |a|^{k/2}}\,
d_\ast a, \qquad 0<\e<\rho<\infty.\ee
 Owing to (\ref{es-wd}), for any test function
$\phi\in\S$, the expression $(I_{\e,\rho}f,\; \phi)$  is finite
when  $f$, $w$, and $\phi$ are
 replaced by $|f|$, $|w|$  and $|\phi|$, respectively. Hence, we
 can  change the order of integration, and (\ref{eq4})   yields
\be \label{aa}
(I_{\e,\rho}f,\; \phi)=(f,\; \dot I_{\e,\rho} \phi).\ee
where $ \dot I_{\e,\rho}$ has the same meaning as in (\ref{Ier}) but with $w$
 replaced by its complex conjugate  $\bar w$. 
Let us show that \be \label{aaa}
(\dot I_{\e,\rho}
\phi)(x)=\F^{-1}[ \, \overline{m_{\e,\rho}} \, \F \phi \, ](x) , \ee where
$m_{\e,\rho}(y)=\tilde  m_{\e,\rho}(y'y)$, \be \tilde
m_{\e,\rho}(r)=\frac{2^{m(k+1)} \pi^{km}}{\sig_{n,m}} \intl_{\{\z\in\Mt: \,\e
r<\z'\z<\rho r\}}\frac{(\tilde\F w)(\z)}{|\z|^{n}_m}\, d\z, \ee $r\in\p$.  Suppose that  the Fourier transform of $w$  has the
form $(\tilde\F w)(\z)=u_0(\z'\z)$ where $u_0$ obeys the symmetry
condition (\ref{sym}). 
By (\ref{eq5}) (with $w$ replaced by
$\overline{w_a}$), we have
 $$
(\dot I_{\e,\rho} \phi)(x)=(2\pi)^{(k-n)m}\!\!\!
\intl_{\vnk}\!\!\!\!d_\ast\xi \!\!\!\intl_{\Mt}\!\!\!\!\!\!
\exp(-\tr(ix'\xi z))k_{\e,\rho}^\a(z)(\F\phi)(\xi z) dz, $$ where
  \bea k_{\e,\rho}^\a(z)&=&\intl_{\e I_m}^{\rho
I_m } \frac{(\tilde \F\bar w)(za^{1/2})}{ |a|^{k/2}}\, d_\ast
a=\intl_{\e I_m}^{\rho I_m }\,  \frac{\overline{ u_0(a^{1/2} ra^{1/2})}}{
|a|^{k/2}}\, d_\ast a, \qquad r=z'z.\nonumber\eea Owing to the
symmetry (\ref{sym}), we have $u_0(a^{1/2} ra^{1/2})=u_0(r^{1/2}
ar^{1/2})$ (without loss of generality, one can assume
$\rank(z)=m$). Then by changing variable $r^{1/2} ar^{1/2}=s$ and
making use of Lemma \ref {l2.3}, we obtain
 \bea
k_{\e,\rho}^\a(z) &=&\intl_{\e I_m}^{\rho I_m } \, \, \frac{  \overline{u_0(r^{1/2}
ar^{1/2})}}{ |a|^{k/2}}\, d_\ast a=|r|^{k/2}\intl_{\e r}^{\rho r } \, 
\frac  {\overline{u_0(s)}}{ |s|^{k/2}}\, d_\ast s \nonumber  \\
&=&\frac{2^m}{\sig_{n-k,m}} \, |r|^{k/2}\intl_{\{\z\in\Mt :\,\e
r<\z'\z<\rho r\}}\frac{(\tilde\F \bar w)(\z)}{|\z|^{n}_m}\, d\z. \nonumber 
 \eea
Replacing $\z$ by $-\z$ and using the equality $(\tilde\F \bar w)(-\z)=
\overline{(\tilde\F w)(\z)}$, we have
\[k_{\e,\rho}^\a(z)=\frac{\sig_{n,m}\, (2\pi)^{-km}}{\sig_{n-k,m}} \, 
 |r|^{k/2} \, \overline{ \tilde
m_{\e,\rho}(r)}. \] Hence
$$
(\dot I_{\e,\rho} \phi)(x)=c\!\!\! \intl_{\vnk}\!\!\!\!d_\ast\xi
\!\!\!\intl_{\Mt}\!\!\!\!\!\! \exp(-\tr(ix'\xi z)) \, |z|_m^k
 \, \overline{\tilde m_{\e,\rho}(z'z)}(\F\phi)(\xi z)\,  dz, $$ where
$c=(2\pi)^{-nm}\sig_{n,m}/\sig_{n-k,m}.$ By Lemma \ref{lSS},
$$ (\dot I_{\e,\rho}
\phi)(x)=(2\pi)^{-nm} \intl_{\Ma} \exp(-\tr(ix'y )) \,\overline{ \tilde
m_{\e,\rho}(y'y)}\, (\F\phi)(y)\; dy ,$$ and (\ref{aaa}) follows.

The rest of the proof is standard.
 Since $f\in L^2$ and  $ m_{\e,\rho}$ is uniformly
  bounded, then  by the Parseval equality,
\be\label {bb} 
 (I_{\e,\rho}f,\; \phi)=(\F^{-1} m_{\e,\rho}\F f,\;
\phi),\qquad \forall \phi\in\S(\Ma).\ee By Lemma \ref{loc}, $
I_{\e,\rho}f$ is a locally integrable function. Since $\F^{-1}
m_{\e,\rho}\F f$ is locally integrable too, then (\ref{bb})
implies a pointwise equality $I_{\e,\rho}f=\F^{-1} m_{\e,\rho}\F
f$, and we have \bea\nonumber \|I_{\e,\rho}f-c_wf\|_2&=&\| \F^{-1}
m_{\e,\rho}\F f-c_wf\|_2\\\nonumber &=& (2\pi)^{-nm/2}\|(
m_{\e,\rho}-c_w)\F f\|_2\to 0\eea as $\e\to 0$, $\rho\to\infty$.
\end{proof}

\subsection{Reproducing formula for  the ridgelet transform}

Given two
 functions $u(z)$ and $v(z)$     on $\Mt$, we set \[
u_a(z)=|a|^{(k-n)/2}u(za^{-1/2}), \qquad
 v_a(z)=|a|^{(k-n)/2}v(za^{-1/2}), \quad a\in\p, \]
 and consider the corresponding
 ridgelet transforms
\bea\label{rtu}(U_a f)(\xi, t)&=& \intl_{\Ma} f(x)\,u_a(t-\xi '
x)\,dx,\\\label{rdtv} (V_a^{\ast}\vp)(x)&=&\intl_{\vnk}\!\!d_\ast
\xi\!\!\intl_{\Mt}\!\!\!\!\!\!\vp(\xi,t)\,v_a(\xi ' x-z)\, dz,\eea
cf. (\ref{rt}), (\ref{rdt}).

\begin{theorem}
Let $u$ and $v$ be integrable radial  functions  on $\Mt$, $1\leq
k\leq n-m$, such that their convolution $u\ast v$ is admissible
(see Definition \ref {adm}). Let \bea \qquad
c_{u,v}&=&\frac{ 2^{m(k+1)} \pi^{km}}{\sig_{n,m}}\intl_{\Mt}\frac{(\tilde\F
u)(\z)(\tilde\F v)(\z)}{|\z|^{n}_m}\, d\z\\\nonumber&=&
 \lim\limits_{A\to 0\atop B\to \infty} \frac {2^{m(k+1)} \pi^{km}}{\sig_{n,m}}
\intl_{\{\z\in\Mt:\,A<z'z<B\}}\frac{(\tilde\F u)(\z)(\tilde\F
v)(\z)}{|\z|^{n}_m}\, dz, \eea ($A,B\in\p$). Then for $f\in
L^1(\Ma)\cap L^2(\Ma)$,
 \be c_{u,v} f=
\intl_{\p } \frac{V_a^\ast U_a f}{ |a|^{k/2}}\, d_\ast a =
\underset{\rho\to\infty }{\lim\limits_{\e \to 0 }^{(L^2)}} \,
\intl_{\e I_m}^{\rho I_m } \frac{V_a^\ast U_a f}{ |a|^{k/2}}\,
d_\ast a . \ee
\end{theorem}
\begin{proof}
Let us show that $V_a^\ast U_a f$ coincides with the dual ridgelet
 transform  $W_a^{\ast}\hat f$ (see (\ref{Wa-r})) generated by the function $w=u\ast
v$. We have \bea (W_a^{\ast}\hat f)(x)&=&\intl_{\vnk}\!\!d_\ast
\xi\!\!\intl_{\Mt}\!\!\!\!\!\!\hat
f(\xi,z)\,dz\!\!\!\!\intl_{\Mt}u_a(t)\,v_a( \xi ' x-z-t)dt\,
\nonumber
\\ &=&\intl_{\vnk}\!\!d_\ast
\xi\!\!\intl_{\Mt}\!\!\!\!\!\!\hat
f(\xi,z)\,dz\!\!\!\!\intl_{\Mt}u_a(\z-z)\,v_a( \xi ' x-\z)d\z
\nonumber
\\ &=&\intl_{\vnk}\!\!d_\ast
\xi\!\!\intl_{\Mt}\!\!\!\!\!\!v_a( \xi ' x-\z)d\z
\!\!\!\!\intl_{\Mt}u_a(\z-z)\hat f(\xi,z)\,dz \nonumber \\
&\stackrel{\rm (\ref{w-r})}{=}&\intl_{\vnk}\!\!d_\ast
\xi\!\!\intl_{\Mt}\!\!\!\!\!\!v_a( \xi ' x-\z)d\z
\!\!\!\!\intl_{\Ma} f(x)\,u_a(\z -\xi ' x)\,dx\nonumber \\
&=&(V_a^\ast U_a f)(x).\nonumber\eea Now the result follows by
Theorem \ref{tinv1}.
\end{proof}

\end{document}